\documentclass[12pt]{article}
\usepackage{epic,eepic,epsf,epsfig}
\usepackage{amsmath}
\usepackage{amssymb}
\usepackage{amsbsy}

\usepackage{amsfonts}
\newtheorem{theorem}{Theorem}[section]%
\newtheorem{lemma}[theorem]{Lemma}%
\newtheorem{cor}[theorem]{Corollary}%

 \def\Omega{\Omega}

\def\f{\noindent}

\newcommand{\qed}{\mbox{\raisebox{0.7ex}{\fbox{}}} \vspace{4truemm}}
\def\demo{\f {\bf Proof.}\hskip10pt}

\begin{document}

\baselineskip 16pt

\title{ \vspace{-1.2cm}
A note on Frobenius quotient for prime-power divisor of the exponent of finite groups
\thanks{\scriptsize This research was supported in part by Shandong Provincial Natural Science Foundation, China (ZR2017MA022)
and NSFC (11761079).
\newline
 \hspace*{0.5cm} \scriptsize $^{\ast\ast}$Corresponding
  author.\newline
       \hspace*{0.5cm} \scriptsize{E-mail addresses:}
       shijt2005@163.com\,(J. Shi),\,lwjytu@qq.com\,(W. Liu).}}

\author{Jiangtao Shi\,$^{\ast\ast}$,\,Wenjing Liu\\
\\
{\small School of Mathematics and Information Sciences, Yantai University, Yantai 264005, China}}

\date{ }

\maketitle \vspace{-.8cm}

\begin{abstract}
Let $G$ be a finite group and $n$ be any prime-power divisor of ${\rm exp}(G)$, the exponent of $G$. Frobenius' theorem indicates that $|\{g\in G\mid g^n=1\}|=f_n\cdot n$ for some positive integer $f_n$.
We call $f_n$ a Frobenius quotient of $G$ for $n$. Let $\mathcal{F}_{pp}(G)=\{f_n\mid n$ is any prime-power divisor of ${\rm exp}(G)$$\}$ and ${\rm mf}_{pp}(G)$ be the maximum Frobenius quotient in $\mathcal{F}_{pp}(G)$. In this paper, we provide a complete classification of finite group $G$ with ${\rm mf}_{pp}(G)\leq q$, where $q$ is the smallest prime divisor of $|G|$.

\medskip

\f {\bf Keywords:} Frobenius quotient; exponent; prime-power divisor; cyclic group; normal\\
{\bf MSC(2010):} 20D10
\end{abstract}

\section{Introduction}\label{s1}

In this paper all groups are assumed to be finite. Let $G$ be a group and $n$ be a divisor of $|G|$. Assume $F_n(G)=\{g\in G\mid\,g^n=1\}$. Frobenius' theorem {\rm\cite{Fro1895}} indicates that $|F_n(G)|=f_n\cdot n$ for some positive integer $f_n$. In {\rm\cite{shi}} $f_n$ is called a Frobenius quotient of $G$ for $n$. Moreover, the set $\mathcal{F}(G)$ that consists of all Frobenius quotients of $G$ is called the Frobenius spectrum of $G$. Let ${\rm mf}(G)$ be the maximum Frobenius quotient of $G$ in $\mathcal{F}(G)$. In {\rm\cite{meng2011}} Meng and Shi obtained a complete classification of group $G$ with ${\rm mf}(G)\leq 2$. Furthermore, Meng, Shi and Chen {\rm\cite{meng2012}} gave a complete classification of group $G$ with ${\rm mf}(G)\leq 3$. {\rm\cite{jiang}} and {\rm\cite{ren}} investigated group $G$ with ${\rm mf}(G)\leq 4$ and ${\rm mf}(G)\leq 5$ respectively. In {\rm\cite{hei20151}} Heineken and Russo characterized group $G$ with $f_n\leq n$ for any divisor $n$ of $|G|$. Moreover, Shi, Meng and Zhang {\rm\cite{shi}} obtained a complete classification of group $G$ with $\mathcal{F}(G)=\{1,\,q\}$, where $q$ is the smallest prime divisor of $|G|$. Let $\mathcal{F}'(G)$ be the set of all Frobenius quotients of $G$ for even divisors of $|G|$ and let ${\rm mfe}(G)$ be the maximum Frobenius quotient of $G$ in $\mathcal{F}'(G)$. Shi, Meng and Zhang {\rm\cite{shi}} showed that if ${\rm mfe}(G)\leq 7$ or ${\rm mfe}(G)\leq 8$ and $A_5$ is not a composition factor of $G$, then $G$ is solvable.

In {\rm\cite{hei20152}} Heineken and Russo denoted by $B(G)={\rm max}\{f_n\mid n\mid{\rm exp}(G)\}$ global breadth of $G$ in the sense of Frobenius, and they {\rm\cite{{hei20152},{hei2019}}} obtained all refined non-abelian groups with $B(G)\leq 11$. Moreover, Meng {\rm\cite{meng2017}} gave a complete characterization of group $G$ with $B(G)=4$.

For a given set $\mathcal{A}$ that consists of some positive integers, let $S_1$ be the set of groups $G$ in which the Frobenius quotient $f_n\in\mathcal{A}$ for any divisor $n$ of $|G|$, $S_2$ be the set of groups $G$ in which the Frobenius quotient $f_n\in\mathcal{A}$ for any divisor $n$ of ${\rm exp}(G)$ and $S_3$ be the set of groups $G$ in which the Frobenius quotient $f_n\in\mathcal{A}$ for any prime-power divisor $n$ of ${\rm exp}(G)$. It is easy to see that $S_1\subseteq S_2\subseteq S_3$.

In this paper, as a further generalization of above literature, our main goal is to investigate group $G$ in which the maximum Frobenius quotient $f_n$ is not greater than a given positive integer for any prime-power divisor $n$ of ${\rm exp}(G)$. Let $\mathcal{F}_{pp}(G)=\{f_n\mid n$ is any prime-power divisor of ${\rm exp}(G)$$\}$ and ${\rm mf}_{pp}(G)$ be the maximum Frobenius quotient in $\mathcal{F}_{pp}(G)$. Our main result is as follows, the proof of which is given in Section~\ref{s3}.

\begin{theorem}\ \ \label{th1} Suppose that $G$ is a group and $q$ is the smallest prime divisor of $|G|$. Let $Q\in{\rm Syl}_q(G)$, then ${\rm mf}_{pp}(G)\leq q$ if and only if $G$ is one of the following groups:

$(1)$ $G=Q\times R$, where $R\geq 1$ is a cyclic $q'$-Hall subgroup of $G$, $Q=Z_{q^m}$ where $m\geq 1$, or $Q=Z_{q^{m-1}}\times Z_q$ where $m\geq 2$, or $Q=Q_8$, or $Q=\langle a, b\mid a^{2^{m-1}}=b^2=1, b^{-1}ab=a^{1+2^{m-2}}\rangle$ where $m\geq 4$, or $Q=\langle a, b\mid a^{q^{m-1}}=b^q=1, b^{-1}ab=a^{1+q^{m-2}}\rangle$ where $q>2$ and $m\geq 3$.

$(2)$ $G=\langle a,\,b\mid\, a^{2^m}=b^3=1,\,a^{-1}ba=b^{-1}\rangle\times T$, where $m\geq 1$, $T$ is cyclic and $(|T|,6)=1$.
\end{theorem}

The following three results are direct corollaries of Theorem~\ref{th1}.

\begin{cor}\ \ \label{c1} Suppose that $G$ is a group and $q$ is the smallest prime divisor of $|G|$. Then ${\rm mf}_{pp}(G)\leq q$ if and only if $\mathcal{F}_{pp}(G)=\{1\}$ or $\{1,\,q\}$.
\end{cor}

\begin{cor}\ \ \label{c2} Suppose that $G$ is a group and $q$ is the smallest prime divisor of $|G|$. If ${\rm mf}_{pp}(G)<q$, then $G$ is cyclic.
\end{cor}

\begin{cor}\ \ \label{c3} Suppose that $G$ is a group and $q$ is the smallest prime divisor of $|G|$. If $\mathcal{F}_{pp}(G)=\{1,\,f\}$, where $f>1$, then $f\geq q$.
\end{cor}

\section{Necessary Lemma}\label{s2}

\begin{lemma}\ \ \label{l1} Suppose that $G$ is a group and $q$ is a prime divisor of $|G|$. Let $Q\in{\rm Syl}_q(G)$ and $|Q|=q^m$, where $m\geq 1$. If $Q$ is cyclic and $G$ has exactly $t$ Sylow $q$-subgroups, then $|F_{|Q|}(G)|=|F_{q^m}(G)|\geq q^m+(t-1)(q^m-q^{m-1})$.
\end{lemma}

\demo Let $Q_1,\,Q_2,\cdots, Q_t$ be $t$ Sylow $q$-subgroups of $G$. It follows that $Q_1\bigcup Q_2\bigcup\cdots\bigcup\\ Q_t\subseteq F_{|Q|}(G)$ and then $|F_{|Q|}(G)|\geq|Q_1\bigcup Q_2\bigcup\cdots\bigcup Q_t|$. It is clear that $|Q_1\bigcup Q_2\bigcup\cdots\bigcup Q_t|\\=|Q_1\bigcup(Q_2\setminus Q_1)\bigcup(Q_3\setminus(Q_1\bigcup Q_2))\bigcup\cdots\bigcup(Q_t\setminus(Q_1\bigcup Q_2\bigcup\cdots\bigcup Q_{t-1}))|=|Q_1|+|Q_2\setminus Q_1|+|Q_3\setminus(Q_1\bigcup Q_2)|+\cdots+|Q_t\setminus(Q_1\bigcup Q_2\bigcup\cdots\bigcup Q_{t-1})|$.

Note that for every $2\leq j\leq t$, one has $|Q_j\setminus(Q_1\bigcup Q_2\bigcup\cdots\bigcup Q_{j-1})|=|Q_j\setminus(Q_j\bigcap(Q_1\bigcup\\ Q_2\bigcup\cdots\bigcup Q_{j-1}))|=|Q_j\setminus((Q_j\bigcap Q_1)\bigcup(Q_j\bigcap Q_2)\bigcup\cdots\bigcup(Q_j\bigcap Q_{j-1}))|$. Since $Q_j$ is a cyclic group of prime-power order, $Q_j$ has exactly one maximal subgroup. Let $Q_{j0}$ be the maximal subgroup of $Q_j$, then $Q_j\bigcap Q_1\leq Q_{j0},\,Q_j\bigcap Q_2\leq Q_{j0},\,\cdots,\,Q_j\bigcap Q_{j-1}\leq Q_{j0}$, which implies that $(Q_j\bigcap Q_1)\bigcup(Q_j\bigcap Q_2)\bigcup\cdots\bigcup(Q_j\bigcap Q_{j-1})\leq Q_{j0}$. Thus $|Q_j\setminus((Q_j\bigcap Q_1)\bigcup(Q_j\bigcap Q_2)\bigcup\cdots\bigcup(Q_j\bigcap Q_{j-1}))|\geq|Q_j\setminus Q_{j0}|=|Q_j|-|Q_{j0}|=q^m-q^{m-1}$.

It follows that $|F_{|Q|}(G)|=|F_{q^m}(G)|\geq|Q_1|+|Q_2\setminus Q_1|+|Q_3\setminus(Q_1\bigcup Q_2)|+\cdots+|Q_t\setminus(Q_1\bigcup Q_2\bigcup\cdots\bigcup Q_{t-1})|\geq q^m+(t-1)(q^m-q^{m-1})$.\hfill\qed

\section{Proof of Theorem~\ref{th1}}\label{s3}

\demo {\bf Part I.} We first prove the necessity part.

{\bf Case 1.} When $G=Q$ is a $q$-group, assume $|G|=q^m$, where $m\geq 1$.

Since ${\rm mf}_{pp}(G)\leq q$, one has ${\rm exp}(G)\geq q^{m-1}$.

If ${\rm exp}(G)=q^m$, then $G$ is a cyclic group.

Next assume ${\rm exp}(G)=q^{m-1}$, then $G$ is a non-cyclic $q$-group having a cyclic maximal subgroup of order $q^{m-1}$.
By {\rm\cite[Theorem 12.5.1]{hall}}, $G$ might be one of the following six classes of groups: (1) $G=\langle a, b\mid a^{q^{m-1}}=b^q=1, [a,b]=1\rangle$, where $m\geq 2$; (2) $G=\langle a, b\mid a^{q^{m-1}}=b^q=1, b^{-1}ab=a^{1+q^{m-2}}\rangle$, where $q>2$ and $m\geq 3$; (3) $G=\langle a, b\mid a^{2^{m-1}}=1, b^2=a^{2^{m-2}}, b^{-1}ab=a^{-1}\rangle$, where $m\geq 3$; (4) $G=\langle a, b\mid a^{2^{m-1}}=b^2=1, b^{-1}ab=a^{-1}\rangle$, where $m\geq 3$; (5) $G=\langle a, b\mid a^{2^{m-1}}=b^2=1, b^{-1}ab=a^{1+2^{m-2}}\rangle$, where $m\geq 4$; (6) $G=\langle a, b\mid a^{2^{m-1}}=b^2=1, b^{-1}ab=a^{-1+2^{m-2}}\rangle$, where $m\geq 4$.

(1) Suppose $G=\langle a, b\mid a^{q^{m-1}}=b^q=1, [a,b]=1\rangle$, where $m\geq 2$. For every $1\leq s\leq m-1$, it is easy to see that $F_{q^s}(G)=\langle a^{q^{m-s-1}}\rangle\times\langle b\rangle$. It follows that $|F_{q^s}(G)|=q\cdot q^s$. Therefore, $f_{q^s}(G)=q$ for every $1\leq s\leq m-1$.

(2) Suppose $G=\langle a, b\mid a^{q^{m-1}}=b^q=1, b^{-1}ab=a^{1+q^{m-2}}\rangle$, where $q>2$ and $m\geq 3$. Note that $Z(G)=\langle a^q\rangle$. It is easy to see that $F_{q^{m-1}}(G)=G$ and $F_{q^s}(G)=\langle a^{q^{m-s-1}}\rangle\times\langle b\rangle$ for every $1\leq s\leq m-2$. Therefore, $f_{q^s}(G)=q$ for every $1\leq s\leq m-1$.

(3) Suppose $G=\langle a, b\mid a^{2^{m-1}}=1, b^2=a^{2^{m-2}}, b^{-1}ab=a^{-1}\rangle$, where $m\geq 3$. If $m\geq 4$, then $F_{2^2}(G)=\langle a^{2^{m-3}}\rangle\bigcup\{a^ib\mid 1\leq i\leq 2^{m-1}\}$. It follows that $|F_{2^2}(G)|=4+2^{m-1}\geq 4+2^3=12>2\cdot 2^2$, this contradicts the hypothesis. If $m=3$, then $G=Q_8$. It is clear that $Q_8$ satisfies the hypothesis.

(4) Suppose $G=\langle a, b\mid a^{2^{m-1}}=b^2=1, b^{-1}ab=a^{-1}\rangle$, where $m\geq 3$. It is easy to see that $F_2(G)=\langle a^{2^{m-2}}\rangle\bigcup\{a^ib\mid 1\leq i\leq 2^{m-1}\}$ and then $|F_2(G)|=2+2^{m-1}\geq 2+2^2=6>2\cdot 2$, this contradicts the hypothesis.

(5) Suppose $G=\langle a, b\mid a^{2^{m-1}}=b^2=1, b^{-1}ab=a^{1+2^{m-2}}\rangle$, where $m\geq 4$. It is easy to see that $o(a^ib)=o(a^i)$ for every $1\leq i\leq 2^{m-1}$. Thus for every $1\leq s\leq m-1$, one has $F_{2^s}(G)=\langle a^{2^{m-1-s}}\rangle\bigcup\langle a^{2^{m-1-s}}b\rangle$ and then $|F_{2^s}(G)|\leq 2\cdot 2^s$, which implies that $f_{2^s}(G)=2$ for every $1\leq s\leq m-1$.

(6) Suppose $G=\langle a, b\mid a^{2^{m-1}}=b^2=1, b^{-1}ab=a^{-1+2^{m-2}}\rangle$, where $m\geq 4$. It is easy to see that $o(a^ib)=2$ if $2\mid i$ and $o(a^ib)=4$ if $2\nmid i$ for every $1\leq i\leq 2^{m-1}$. Thus $F_2(G)=\langle a^{2^{m-2}}\rangle\bigcup\{a^ib\mid i=2r,\,1\leq r\leq 2^{m-2}\}$ and then $|F_2(G)|=2+2^{m-2}\geq 2+4=6>2\cdot 2$, this contradicts the hypothesis.

{\bf Case 2.} When $G$ is not a $q$-group, assume $|G|=q^m{q_2}^{m_2}\cdots{q_k}^{m_k}$, where $q<q_2<\cdots<q_k$, $k\geq 2$, $m\geq 1$ and $m_i\geq 1$ for every $2\leq i\leq k$.

For every $q_i$ where $2\leq i\leq k$, let $Q_i\in{\rm Syl}_{q_i}(G)$, we claim that $Q_i$ is cyclic.

Otherwise, if $Q_i$ is non-cyclic, take $n={\rm exp}(Q_i)<|Q_i|$. Then $|F_n(G)|\geq|Q_i|\geq q_i\cdot n>q\cdot n$, a contradiction. Therefore, $Q_i$ is cyclic for every $2\leq i\leq k$.

Next we claim $Q_i\trianglelefteq G$ for every $2\leq i\leq k$.

Otherwise, assume $Q_i\ntrianglelefteq G$ for some $2\leq i\leq k$. Take $n={\rm exp}(Q_i)=|Q_i|={q_i}^{m_i}$, where $m_i\geq 1$. Note that $G$ has at least $q_i+1$ Sylow $q_i$-subgroups by Sylow's theorem. By Lemma~\ref{l1}, one has $|F_n(G)|\geq {q_i}^{m_i}+(q_i+1-1)({q_i}^{m_i}-{q_i}^{m_i-1})={q_i}^{m_i}+q_i({q_i}^{m_i}-{q_i}^{m_i-1})=q_i\cdot {q_i}^{m_i}>q\cdot {q_i}^{m_i}$, this contradicts the hypothesis. Therefore, $Q_i\trianglelefteq G$ for every $2\leq i\leq k$.

Let $Q\in{\rm Syl}_q(G)$,where $|Q|=q^m$. Since $q\cdot{{\rm exp}(Q)}\geq|F_{{\rm exp}(Q)}(G)|\geq |F_{{\rm exp}(Q)}(Q)|\geq|Q|$, one has ${\rm exp}(Q)=q^m$ or $q^{m-1}$.

For the case when ${\rm exp}(Q)=q^{m-1}$, we claim $Q\trianglelefteq G$. Otherwise, assume $Q\ntrianglelefteq G$. Then $G$ has more than one Sylow $q$-subgroups, which implies that $|F_{{\rm exp}(Q)}(G))|>|Q|=q\cdot {\rm exp}(Q)$, this contradicts the hypothesis. Therefore, $Q\trianglelefteq G$.

Combine the analyses in Case 1, one has $G=Q\times Q_2\times\cdots\times Q_k$, where $Q=Z_{q^{m-1}}\times Z_q$ for $m\geq 2$, or $Q=Q_8$, or $Q=\langle a, b\mid a^{2^{m-1}}=b^2=1, b^{-1}ab=a^{1+2^{m-2}}\rangle$ for $m\geq 4$, or $Q=\langle a, b\mid a^{q^{m-1}}=b^q=1, b^{-1}ab=a^{1+q^{m-2}}\rangle$ for $q>2$ and $m\geq 3$, $Q_i$ is cyclic for every $2\leq i\leq k$.

For another case when ${\rm exp}(Q)=q^m$, one has that $Q$ is cyclic.

$(i)$ Suppose $Q \trianglelefteq G$, then $G$ is a cyclic group.

$(ii)$ Suppose $Q\ntrianglelefteq G$. By Sylow's theorem, assume that $G$ has exactly $hq+1$ Sylow $q$-subgroups, where $h\geq 1$. By Lemma~\ref{l1}, one has $|F_{|Q|}(G)|\geq q^m+(hq+1-1)(q^m-q^{m-1})=(hq-h+1)\cdot q^m$. Note that $hq-h+1\geq q$ and $hq-h+1=q$ if and only if $h=1$. By the hypothesis, $G$ has exactly $q+1$ Sylow $q$-subgroups. Moreover, one has $q+1\mid|G|$ by Sylow's theorem. Since $q$ is the smallest prime divisor of $|G|$, $q+1$ must be a prime divisor of $|G|$. It follows that $q=2$ and $q+1=3$.

Let $Q_2$ be a Sylow $3$-subgroup of $G$. It is easy to see that $QQ_i=Q\times Q_i$ for every $3\leq i\leq k$ and $QQ_2=Q\ltimes Q_2$, where $|G:N_G(Q)|=|QQ_2:N_{QQ_2}(Q)|=3$.

In the following we will show that $Q\ltimes Q_2$ is a minimal non-nilpotent group.

Note that $QQ_2$ is a metacyclic group and then $QQ_2$ is supersolvable. It follows that every maximal subgroup of $QQ_2$ has index a prime in $QQ_2$.

Let $M$ be a maximal subgroup of $QQ_2$ such that $|QQ_2:M|=3$, then $Q\leq M$. One has $M=M\cap QQ_2=Q(M\cap Q_2)=Q\ltimes(M\cap Q_2)$. It follows that $M\cap Q_2$ is maximal in $Q_2$. Since $Q_2$ is cyclic, $M\cap Q_2$ is the unique maximal subgroup of $Q_2$. Then $M$ is the unique maximal subgroup of $QQ_2$ that has index 3. It follows that $N_{QQ_2}(Q)=M$ and then $M=Q\times(M\cap Q_2)$ is nilpotent.

Let $N$ be a maximal subgroup of $QQ_2$ that has index 2 in $QQ_2$, then $Q_2\leq N$. One has $N=N\cap QQ_2=(N\cap Q)\ltimes Q_2$, where $N\cap Q$ is the unique maximal subgroup of $Q$ since $Q$ is a cyclic group of prime-power order. Let $Q$, $Q_{01}$ and $Q_{02}$ be 3 Sylow 2-subgroups of $G$. Since $2\cdot 2^m\geq|F_{2^m}(G)|=|Q\bigcup(Q_{01}\setminus Q)\bigcup(Q_{02}\setminus(Q\bigcup Q_{01}))|=|Q\bigcup(Q_{01}\setminus(Q\cap Q_{01}))\bigcup(Q_{02}\setminus((Q_{02}\cap Q)\bigcup(Q_{02}\cap Q_{01}))|=|Q|+|Q_{01}\setminus(Q\cap Q_{01})|+|Q_{02}\setminus((Q_{02}\cap Q_{01})\bigcup(Q_{02}\cap Q_{01}))|\geq |Q|+(|Q_{01}|-|N_{01}|)+(|Q_{02}|-|N_{02}|)=2^m+(2^m-2^{m-1})+(2^m-2^{m-1})=2\cdot 2^m=2\cdot|Q|$, where $N_{01}$ is the unique maximal subgroup of $Q_{01}$ and $N_{02}$ is the unique maximal subgroup of $Q_{02}$. It follows that $Q_{01}\cap Q=N_{01}$ and $(Q_{02}\cap Q)\bigcup(Q_{02}\cap Q_{01})=N_{02}$. First, one has $Q_{01}\cap Q=N\cap Q$ since both $N\cap Q$ and $Q_{01}\cap Q$ are maximal subgroups of $Q$.
Second, since every group cannot be written as an union of two proper subgroups, one has either $Q_{02}\cap Q=N_{02}$ or $Q_{02}\cap Q_{01}=N_{02}$. If $Q_{02}\cap Q=N_{02}$, then $Q_{02}\cap Q=N\cap Q$ since both $N\cap Q$ and $Q_{02}\cap Q$ are maximal subgroups of $Q$. If $Q_{02}\cap Q_{01}=N_{02}$, then $Q_{02}\cap Q_{01}=Q_{01}\cap Q$ since both $Q_{02}\cap Q_{01}$ and $Q_{01}\cap Q$ are maximal subgroups of $Q_{01}$. We also have $Q_{02}\cap Q_{01}=N\cap Q$. Thus $Q_G=Q\cap Q_{01}\cap Q_{02}=(Q\cap Q_{01})\cap(Q\cap Q_{02})\cap(Q_{01}\cap Q_{02})=N\cap Q$. It follows that $N\cap Q\trianglelefteq G$ and then $N=(N\cap Q)\times Q_2$ is nilpotent.

By above arguments, one has that $Q\ltimes Q_2$ is a minimal non-nilpotent group.

By {\rm\cite[Chapter 3, Theorem 5.2]{huppert}}, one has ${\rm exp}(Q_2)=3$. Since $Q_2$ is cyclic, one has $|Q_2|=3$. Let $Q=\langle a\rangle,\,Q_2=\langle b\rangle$, one has $Q\ltimes Q_2=\langle a\rangle\ltimes\langle b\rangle=\langle a,\,b\mid\, a^{-1}ba=b^{-1}\rangle$. It follows that  $G=\langle a,\,b\mid\, a^{2^m}=b^3=1,\,a^{-1}ba=b^{-1}\rangle\times(Q_3\times\cdots\times Q_k)$, where $Q_3,\cdots,\,Q_k$ are cyclic.

{\bf Part II.} According to the above analyses, the sufficiency part is evident.\hfill\qed

\end{document}